\documentclass[12pt,draftcls,onecolumn] {IEEEtran} % Single Column for review

\usepackage[T1]{fontenc}
\usepackage[latin9]{inputenc}
\usepackage{geometry}
\usepackage{units}
\usepackage{mathrsfs}
\usepackage{amsmath}
\usepackage{amsthm}
\usepackage{amssymb}
\usepackage{graphicx}
\usepackage{caption}
\usepackage{subcaption}
\usepackage{color}

\makeatletter
\theoremstyle{plain}
\newtheorem{thm}{\protect\theoremname}
\theoremstyle{plain}
\newtheorem{lem}[thm]{\protect\lemmaname}
\theoremstyle{remark}
\newtheorem{rem}[thm]{\protect\remarkname}
\theoremstyle{plain}
\newtheorem{prop}[thm]{\protect\propositionname}

\makeatother

\usepackage{babel}
\providecommand{\lemmaname}{Lemma}
\providecommand{\propositionname}{Proposition}
\providecommand{\remarkname}{Remark}
\providecommand{\theoremname}{Theorem}

\begin{document}
\pdfminorversion=4
\title{Arbitrary Order Fixed-Time Differentiators}
\date{\today}
\author{Jaime A. Moreno,~\IEEEmembership{Member~IEEE} 
\thanks{J.~A. Moreno is with Instituto de Ingenier\'{\i}a, Universidad Nacional Aut\'{o}noma de M\'{e}xico, 04510 Coyoac\'{a}n, Mexico City, Mexico. E-mail: JMorenoP@ii.unam.mx}% <-this % stops a space
\thanks{Financial support PAPIIT-UNAM, project IN110719.} 
}
\maketitle
\begin{abstract}
Differentiation is an important task in control, observation and fault detection. Levant's differentiator is unique, since it is able to estimate exactly and robustly the derivatives of a signal with a bounded high-order derivative. However, the convergence time, although finite, grows unboundedly with the norm of the initial differentiation error, making it uncertain when the estimated derivative is exact. In this paper we propose an extension of Levant's differentiator so that the worst case convergence time can be arbitrarily assigned independently of the initial condition, i.e. the estimation converges in \emph{Fixed-Time}. We propose also a family of continuous differentiators and provide a unified Lyapunov framework for analysis and design. 
\end{abstract}

\section{Introduction}
\label{Sec:Intro}
Given a (Lebesgue-measurable) signal $f\left(t\right)$
defined on $[0,\infty)$ the objective of a differentiator is to estimate
as close as possible some of its time derivatives. Usually, signal
$f\left(t\right)$ is composed of the base signal $f_{0}\left(t\right)$,
which we want to differentiate and is assumed to be $n$-times differentiable,
and a noise signal $\nu(t)$, that we will assume to be uniformly
bounded, i.e. $f(t)=f_{0}(t)+\nu(t)$. 

In order to estimate the derivatives $f_{0}^{\left(i\right)}\left(t\right)=\frac{d^{i}}{dt^{i}}f_{0}\left(t\right)$,
for $i=1,\cdots,\,n-1$, we propose the following nonlinear family
of differentiators ($i=1,\cdots,\,n-1$)
\begin{equation} \label{eq:LevDif_n}
\begin{split}
\dot{x}_{i} & = -k_{i}\phi_{i}\left(x_{1}-f\right)+x_{i+1}\,,\\
\dot{x}_{n} & =  -k_{n}\phi_{n}\left(x_{1}-f\right)\,,
\end{split} 
\end{equation}
where the nonlinear output injection terms, given by 
\begin{equation}
\phi_{i}\left(z\right)=\varphi_{i}\circ\cdots\varphi_{2}\circ\varphi_{1}\left(z\right)\,,\label{eq:DefPhi_i}
\end{equation}
are the composition of the monotonic growing functions $\varphi_{i}:\mathbb{R}\rightarrow\mathbb{R}$ (note that $\lfloor z\rceil^{p}=|z|^{p}\mbox{sign}(z)$)
\begin{equation}
\varphi_{i}\left(s\right)=\kappa_{i}\left\lceil s\right\rfloor ^{\frac{r_{0,\,i+1}}{r_{0,\,i}}}+\theta_{i}\left\lceil s\right\rfloor ^{\frac{r_{\infty,\,i+1}}{r_{\infty,\,i}}}\,.\label{eq:DefVarphi_i}
\end{equation}
$\varphi_{i}$ is a sum of two (signed) power functions, with powers
selected as $r_{0,\,n} =r_{\infty,\,n} =1$, and for $i=1,\cdots,\,n+1$
\begin{equation} \label{eq:Def_r}
\begin{split}
r_{0,\,i} & =  r_{0,\,i+1}-d_{0} = 1-\left(n-i\right)d_{0}\,, \\
r_{\infty,\,i} & =  r_{\infty,\,i+1}-d_{\infty} = 1-\left(n-i\right)d_{\infty}\,, 
\end{split}
\end{equation}
which are completely defined by two parameters $-1\leq d_{0}\leq d_{\infty}<\frac{1}{n-1}$.
With this selection the powers in (\ref{eq:DefVarphi_i}) satisfy
$\frac{r_{0,\,i+1}}{r_{0,\,i}}\leq\frac{r_{\infty,\,i+1}}{r_{\infty,\,i}}$,
so that the first term in $\varphi_{i}\left(s\right)$ is dominating
for small values of $s$, while the second 
is dominating for large values of $s$. This domination effect is
naturally extended to the injection terms $\phi_{i}$ in (\ref{eq:DefPhi_i}).
The (internal) gains $\kappa_{i}>0$ and $\theta_{i}>0$ can be selected
as arbitrarily positive values, and correspond to the desired weighting
of each of the terms in $\varphi_{i}$ (and therefore in $\phi_{i}$).
One possible and simple selection is $\kappa_{i}=\mu$ and $\theta_{i}=1-\mu$
for $i=1,\cdots,\,n$, with $0<\mu<1$ giving the weight of the low-power
and the high-power terms. %
Note that, since for $d_{0}=-1$ system (\ref{eq:LevDif_n})
has a discontinuous right hand side, their solutions are understood
in the sense of Filippov \cite{Fil88}.

Some well-known differentiators in the literature are homogeneous.
For example, the High-Gain observer used as differentiator in \cite{VasKha08,PraKha13}
(see also \cite{KhaPra14}), being linear, is homogeneous of degree
zero. The classical robust and exact differentiator proposed by Levant
\cite{Lev98,Lev03,Lev05} (see also \cite{CruMor19}),
has discontinuous injection terms and is also homogeneous. A family
of homogeneous differentiators, including the previous ones, has been
also proposed recently \cite{CruMor16,SanCru18,CruMor19} for non positive
homogeneity degrees, and in \cite{JbaLev20} for arbitrary degrees.

Differentiator \eqref{eq:LevDif_n} is not homogeneous, but
it is homogeneous in the bi-limit \cite{AndPra08} (bl-homogeneous
for short), that is, near to the origin it is approximated by a homogeneous
system of degree $d_{0}$ and far from the origin it is approximated
by a homogeneous system of degree $d_{\infty}$. Although the scaling
properties of the homogeneous systems are lost, the design of bl-homogeneous
differentiators is more flexible, since the properties near the origin
and far from it can be assigned independently. In particular, by selecting
$d_{0}=d_{\infty}=d$ the differentiator \eqref{eq:LevDif_n}
becomes homogeneous. For $d=0$ one obtains the High-Gain
differentiator, for $d=-1$ Levant's robust and exact
differentiator is recovered and for other values of $d$
the family of differentiators in \cite{CruMor16,SanCru18,CruMor19,JbaLev20}
is attained. Note that for $d<0$ (resp. $d=0$) the estimation converges
in finite-time (resp. exponentially). For $d>0$ the convergence is asymptotic, but
it attains any neighborhood of zero in a time which is uniform in
the initial conditions \cite{AndPra08}. 

Of particular interest for
a differentiator is a property that is only achieved when $d_{0}=d_{\infty}=-1$ \cite{Lev98,Lev03,Lev05}.
In that case $\phi_{n}$ is discontinuous and it induces a Higher-Order
Sliding-Mode at the origin, allowing the estimation to converge (in
the absence of noise) exactly, robustly and in finite-time to the
actual values of the signal derivatives, when the $n$-th derivative
of the signal is bounded by a non zero constant, i.e. $\left|f_{0}^{\left(n\right)}\left(t\right)\right|\leq\Delta$.
For all other values of $d_{0}=d_{\infty}>-1$, convergence is only
achieved if $\Delta=0$. 

One disadvantage of homogeneous (including Levant's exact) differentiators
with negative homogeneity degree, is that the convergence time, although
finite, grows unboundedly (and faster than linearly) with the size
of the initial estimation error. One of the nice features of the bl-homogeneous
design in general \cite{AndPra08}, and of the proposed differentiator
\eqref{eq:LevDif_n} in particular, is that assigning a positive
homogeneity degree to the $\infty$-limit approximation $d_{\infty}>0$
and a negative homogeneity degree to the $0$-limit approximation
$d_{0}<0$, it is possible to counteract this effect: Convergence
of the estimation will be achieved in \emph{Fixed-Time (FxT)} \cite{LopPol18},
that is, the estimation error converges globally, in finite-time and
the settling-time function is globally bounded by a positive constant
$\overline{T}$, independent of the initial estimation error. This
is an important feature, since the differentiator can be designed such that
we are sure that after an arbitrarily assigned time $\overline{T}$ the estimation is
correct no matter what the initial conditions are. Moreover, if $d_{0}=-1$
exact and robust estimation is obtained for all signals having bounded
Lipschitz constant $\left|f_{0}^{\left(n\right)}\left(t\right)\right|\leq\Delta$,
and not only for time polynomial signals, for which $f_{0}^{\left(n\right)}\left(t\right)\equiv0$.
For the first order differentiator (i.e. $n=2$) this property has
been obtained in \cite{CruMor10,CruMor11,Mor13,FraAng12}, using
quadratic-like Lyapunov functions \cite{Mor09,Mor11}. This approach
for the discontinuous first-order differentiator has been extended
and refined in \cite{SeeHai20}, where a detailed gain scaling has
been developed and a tight convergence time estimation has been obtained,
using the results of \cite{SeeHor18}. For differentiators of arbitrary
order this can be achieved by using a switching strategy between two
homogeneous differentiators of positive and negative degrees, as it
is proposed in \cite{AngMor13}. In \cite{LopPol18} also a switching
strategy between homogeneous differentiators with restricted
degrees is presented.

This work can be seen as an extension to an arbitrary
order of the smooth strategy of combining two homogeneous differentiators
proposed in \cite{CruMor10,CruMor11,Mor13}, and in the recent work
\cite{SeeHai20}. Our construction extends to the discontinuous case
the recursive observer design developed for continuous homogeneous
observers in \cite{YanLin04,QiaLin06} and highly improved in \cite{AndPra08,AndPra09} for \emph{continuous} bl-homogeneous observers.
Although many combinations in the selection of $d_{0}\leq d_{\infty}$
are possible, we are particularly interested in the cases $-1\leq d_{0}\leq0\leq d_{\infty}<\frac{1}{n-1}$,
and especially when $d_{0}=-1$.

%The main contribution of this work is to extend the Lyapunov method to the discontinuous arbitrary order Levant\textquoteright s differentiator for the first time.

%The rest of the paper is organized as follows. 
In Section \ref{sec:Preliminaries}, some necessary concepts on bl-homogeneous functions and systems are briefly recalled. Section \ref{sec:Differentiator} presents the main properties of the proposed differentiator. Section \ref{sec:LyapunovProof} contains all the proofs. A simulation example is presented in Section \ref{sec:Examples}. In Section \ref{sec:Conclusions} we draw some conclusions.

%%%%%%%%%%%%%%%%%%%%%%%%%%%%%%%%%%%%%%%%%%%%
\section{Preliminaries}
\label{sec:Preliminaries}
Our notation is fairly standard. We recall briefly some definitions of homogeneity and homogeneity in the bi-limit. However, for precise definitions and properties, we refer the reader to \cite{BacRos05,Bhat2005,BerEfi14} for homogeneity of continuous or discontinuous systems, and to \cite{AndPra08,CruMor20b} for homogeneity in the bi-limit of continuous or discontinuous systems, respectively.
  
For a vector $x\in\mathbb{R}^{n}$, all real values $\epsilon>0$, and $n$ positive real numbers $r_{i}>0$ the dilation operator is defined as $\Delta_{\epsilon}^{\mathbf{r}}x=\left[\epsilon^{r_{1}}x_{1},...,\,\epsilon^{r_{n}}x_{n}\right]^{\top}$. Constants $r_{i}>0$ are the weights of the coordinates $x_{i}$, and $\textbf{r}:=\left[r_{1},...,\,r_{n}\right]$ is the vector of weights. %
A function $V:\mathbb{R}^{m}\mapsto\mathbb{R}^{n}$ (resp. a vector field $f:\mathbb{R}^{n}\mapsto\mathbb{R}^{n}$) is said to be ${\bf r}$-homogeneous of degree $l \in\mathbb{R}$, or $({\bf r},l)$-homogeneous for short, if for all $\epsilon>0$ and for all $x\in\mathbb{R}^{m}\setminus\{0\}$ the equality $V\left(\Delta_{\epsilon}^{{\bf r}}x\right)=\epsilon^{l}V\left(x\right)$ (resp., $f(\Delta_{\epsilon}^{{\bf r}}x)=\epsilon^{l}\Delta_{\epsilon}^{{\bf r}}f(x)$) holds.

A function $\varphi:\mathbb{R}^{n}\mapsto\mathbb{R}$ is said to be homogeneous in the $0$-limit with associated triple $\left(\mathbf{r}_{0},\,l_{0},\,\varphi_{0}\right)$, if it is approximated near $x=0$ by the $({\bf r_{0}},l_{0})$-homogeneous function $\varphi_{0}$. It is said to be homogeneous in the $\infty$-limit with associated triple $\left(\mathbf{r}_{\infty},\,l_{\infty},\,\varphi_{\infty}\right)$, if it is approximated near $x=\infty$ by the $({\bf r_{\infty}},l_{\infty})$-homogeneous function $\varphi_{\infty}$. Similar definitions apply for vector fields and set-valued vector fields. Finally, a function $\varphi:\mathbb{R}^{n}\mapsto\mathbb{R}$ (or a vector field $f:\mathbb{R}^{n}\rightarrow\mathbb{R}^{n}$, or set-valued vector field $F:\mathbb{R}^{n}\rightrightarrows\mathbb{R}^{n}$) is said to be homogeneous in the bi-limit, or \emph{bl-homogeneous} for short, if it is homogeneous in the $0$-limit and homogeneous in the $\infty$-limit. 

%%%%%%%%%%%%%%%%%%%%%%%%%%%%%%%%%%%%%%%%%%%%
\section{Properties of the differentiator}

\label{sec:Differentiator}
The main result of this work states that the differentiator \eqref{eq:LevDif_n},
in the absence of noise, is able to estimate asymptotically the first
$n-1$ derivatives of the signal $f_{0}\left(t\right)$. Let $\mathscr{S}_{0}^{n} \triangleq \left\{ f^{\left(n\right)}\left(t\right)\equiv0\right\} $ represent the class of polynomial signals, while $\mathscr{S}_{\Delta}^{n} \triangleq \left\{ \left|f^{\left(n\right)}\left(t\right)\right|\leq\Delta\right\} $ corresponds to the the class of $n$-Lipschitz
signals.
\begin{thm}
\label{thm:Differentiator}Assume that the signal $f\left(t\right)$
satisfy the stated conditions. Select $-1\leq d_{0}\leq d_{\infty}<\frac{1}{n-1}$
and choose arbitrary positive (internal) gains $\kappa_{i}>0$ and
$\theta_{i}>0$, for $i=1,\cdots,\,n$. Assume further that $\left|f_{0}^{\left(n\right)}\left(t\right)\right|\leq\Delta$
for some non negative Lipschitz constant $\Delta>0$ if $d_{0}=-1$
or $\Delta=0$ if $d_{0}>0$. Under these conditions and in the absence
of noise ($\nu\left(t\right)\equiv0$), there exist appropriate gains
$k_{i}>0$, for $i=1,\cdots,\,n$, such that the solutions of the
bl-homogeneous differentiator \eqref{eq:LevDif_n} converge globally
and asymptotically to the derivatives of the signal, i.e. $x_{i}\left(t\right)\rightarrow f_{0}^{\left(i-1\right)}\left(t\right)$
as $t\rightarrow\infty$. In particular, they converge in \emph{Fixed-Time}, i.e. $x_{i}\left(t\right)\rightarrow f_{0}^{\left(i-1\right)}\left(t\right)$
as $t\rightarrow \bar{T}$, for $i=1,\cdots,\,n$, if either \\
(a) $-1<d_{0}<0<d_{\infty}<\frac{1}{n-1}$ and $f \left(t\right) \in \mathscr{S}_{0}^{n}$, or \\
(b) $-1=d_{0}<0<d_{\infty}<\frac{1}{n-1}$ and $f \left(t\right) \in \mathscr{S}_{\Delta}^{n}$. \qed
\end{thm}
All proofs are given in Section
\ref{sec:LyapunovProof}. The distinguishing feature of the differentiator
\eqref{eq:LevDif_n}, compared to their homogeneous counterparts,
is that it converges within a Fixed-Time when $d_{0}<0<d_{\infty}$.
For the discontinuous differentiator ($d_{0}=-1$) this is accomplished
for a much larger class of signals, since $\mathscr{S}_{0}^{n}\subset\mathscr{S}_{\Delta}^{n}$,
and $\mathscr{S}_{\Delta}^{n}$ is much larger than $\mathscr{S}_{0}^{n}$. 
\begin{rem}
The selection of the function $\varphi_{i}$ in (\ref{eq:DefVarphi_i})
is dictated by the simplicity and concreteness of the presentation.
However, a rather large family of functions can be selected if they satisfy similar appropriate conditions.
\end{rem}

\subsection{Differentiation error dynamics and Lyapunov function}

Defining the differentiation error as $e_{i}\triangleq x_{i}-f_{0}^{\left(i-1\right)}$,
their dynamics satisfy ($i=1,\cdots,\,n-1$)
\begin{equation} \label{eq:Dif_error}
\begin{split}
\dot{e}_{i} & =  -k_{i}\phi_{i}\left(e_{1}-\nu\right)+e_{i+1}\,, \\
\dot{e}_{n} & =  -k_{n}\phi_{n}\left(e_{1}-\nu\right)+\delta\left(t\right)\,,
\end{split} 
\end{equation}
where $\delta\left(t\right)=-f_{0}^{\left(n\right)}\left(t\right)$.
For the variables $z_{1}=\frac{e_{1}}{1}$, $z_{i}=\frac{e_{i}}{k_{i-1}}$ for $i=2,\cdots,\,n$, %
%\[
%z_{1}=\frac{e_{1}}{1},\,z_{2}=\frac{e_{2}}{k_{1}},\cdots,\,z_{i}=\frac{e_{i}}{k_{i-1}},\,\cdots,\,z_{n}=\frac{e_{n}}{k_{n-1}},
%\]
the dynamics of (\ref{eq:Dif_error}) becomes
\begin{equation} \label{eq:ErrorDyn}
\begin{split}
\dot{z}_{i} & =  -\tilde{k}_{i}\left(\phi_{i}\left(z_{1}-\nu\right)-z_{i+1}\right)\,, \\
\dot{z}_{n} & =  -\tilde{k}_{n}\left(\phi_{n}\left(z_{1}-\nu\right)-\bar{\delta}\left(t\right)\right)\,,
\end{split} 
\end{equation}
where for $i=1,\cdots,\,n,$ 
\[
\tilde{k}_{i}=\frac{k_{i}}{k_{i-1}},\,k_{0}=1\,,\,\bar{\delta}\left(t\right)=-\frac{f^{\left(n\right)}\left(t\right)}{k_{n}}\,.
\]

For the convergence proof we will use a (smooth) bl-homogeneous Lyapunov
Function $V$. To define it we select for $n\geq2$ two positive real numbers
$p_{0}$ and $p_{\infty}$, corresponding to the homogeneity degrees
of the $0$-limit and the $\infty$-limit approximations of $V$, such that
\begin{equation} \label{eq:Cond_p1}
\begin{split}
p_{0} & \geq\max_{i\in\left\{ 1,\cdots,\,n\right\} }\left\{ \frac{r_{0,\,i}}{r_{\infty,\,i}}\left(2r_{\infty,\,i}+d_{\infty}\right)\right\} \,, \\\,p_{\infty} & \geq 2\max_{i\in\left\{ 1,\cdots,\,n\right\} }\left\{ r_{\infty,\,i}\right\} +d_{\infty}\,,
\end{split}
\end{equation}
\begin{equation}
\frac{p_{0}}{r_{0,\,i}}<\frac{p_{\infty}}{r_{\infty,\,i}}\,.\label{eq:Cond_p2}
\end{equation}
For $i=1,\,\cdots,\,n$ choose arbitrary positive real
numbers $\beta_{0,\,i}>0$, $\beta_{\infty,\,i}>0$ and define the
functions
\begin{gather} 
Z_{i}\left(z_{i},\,z_{i+1}\right)  = \sum_{j=\{0,\,\infty \}} \qquad \qquad \qquad \qquad \label{eq:DefZi} \\
\beta_{j,i} \left[ \frac{r_{j,i}}{p_{j}}\left|z_{i}\right|^{\frac{p_{j}}{r_{j,i}}} - z_{i}\left\lceil \xi \right\rfloor ^{\frac{p_{j}-r_{j,i}}{r_{j,i}}}+ \frac{p_{j}-r_{j,i}}{p_{j}} \left| \xi \right|^{\frac{p_{j}}{r_{j,i}}} \right]  \nonumber
\end{gather}
where $\xi = \varphi_{i}^{-1}\left(z_{i+1}\right)$, and for $i=1,\,\cdots,\,n-1$, $\varphi_{i}^{-1}$
is the inverse function of $\varphi_{i}$ (\ref{eq:DefVarphi_i}). %
For $i=n$ take $\xi = z_{n+1}\equiv0$, i.e. $Z_{n}\left(z_{n}\right)=\beta_{0,\,n}\frac{1}{p_{0}}\left|z_{n}\right|^{p_{0}}+\beta_{\infty,\,n}\frac{1}{p_{\infty}}\left|z_{n}\right|^{p_{\infty}}$.
The Lyapunov Function candidate is then defined as
\begin{eqnarray}
V\left(z\right) & = & \sum_{j=1}^{n-1}Z_{j}\left(z_{j},\,z_{j+1}\right)+Z_{n}\left(z_{n}\right)\,.\label{eq:DefLF}
\end{eqnarray}

\begin{prop}
\label{prop:LyapunovF}Let the hypothesis of Theorem \ref{thm:Differentiator}
be satisfied, and select $p_{0}$ and $p_{\infty}$ such that (\ref{eq:Cond_p1}) and (\ref{eq:Cond_p2}) are fulfilled.
Under these conditions and in the absence of noise, there exist gains
$k_{i}>0$, for $i=1,\cdots,\,n$, such that $V\left(z\right)$ in
(\ref{eq:DefLF}) is a $\mathcal{C}^{1}$, bl-homogeneous Lyapunov
function for the estimation error dynamics \eqref{eq:ErrorDyn} for
any selection $-1\leq d_{0}\leq d_{\infty}<\frac{1}{n-1}$ and if
$\Delta=0$ in case $d_{0}\neq-1$. Moreover, $V$ satisfies the differential
inequality (\ref{eq:LyapIneq}) for some positive constants $\eta_{0}$,
$\eta_{\infty}$
\begin{equation}
\dot{V}\left(z\right)\leq-\eta_{0}V^{\frac{p_{0}+d_{0}}{p_{0}}}\left(z\right)-\eta_{\infty}V^{\frac{p_{\infty}+d_{\infty}}{p_{\infty}}}\left(z\right)\,.\label{eq:LyapIneq}
\end{equation}
Thus, $z=0$ is a Globally Asymptotically Stable equilibrium point
of (\ref{eq:ErrorDyn}). In particular, %
if $d_{0}<0<d_{\infty}$ then $z=0$ is \emph{Fixed-Time Stable
(FxTS)} \cite{Polyakov2012}, that is, it is globally FTS and the
settling-time function $T\left(z_{0}\right)$ is globally bounded
by a positive constant $\overline{T}$, independent of $z_{0}$, i.e.,
$\exists\overline{T}\in\mathbb{R}_{>0}$ such that $\forall z_{0}\in\mathbb{R}^{n}$, 
$T\left(z_{0}\right)\leq\overline{T}$. %
 $T\left(\cdot\right)$ is continuous at zero and locally bounded.
Moreover, the Fixed-Time $\overline{T}$ can be estimated from \eqref{eq:LyapIneq} as
\begin{equation}
\bar{T} \leq \frac{p_{0}}{d_{0}\eta_{\infty}} \left(\frac{p_{\infty}d_{0}}{p_{0}d_{\infty}}-1\right) \left( \frac{\eta_{0}}{\eta_{\infty}}\right)^{\frac{1}{\left( \frac{p_{\infty}d_{0}}{p_{0}d_{\infty}} -1\right)}} \,.\label{eq:FxTEstimation}
\end{equation}
\qed
\end{prop}
Theorem \ref{thm:Differentiator} is in fact a consequence of this proposition.

\subsection{Gain Calculation}

Stabilizing gains $k_{i}>0$, $i=1,\cdots,\,n$, for the differentiator \eqref{eq:ErrorDyn} can be calculated using $V\left(z\right)$.
\begin{prop}
\label{prop:GainCalculation} Let the hypothesis of Proposition \ref{prop:LyapunovF} be satisfied. A sequence of stabilizing gains $k_{i}>0$, for $i=1,\cdots,\,n$,
can be calculated backwards as follows:

(a) Select $k_{n}\kappa_{n}>\Delta$ and $\tilde{k}_{n}>0$.

(b) For $i=n-1,\,n-2,\cdots,\,1$ select 
\begin{eqnarray*}
\tilde{k}_{i} & > & \omega_{i}\left(\tilde{k}_{i-1},\cdots,\,\tilde{k}_{n}\right)\,.
\end{eqnarray*}
Functions $\omega_{i}$ are given by \eqref{Def:Omegai}, in Section
\ref{subsec:Gain-calculation}, and are obtained from $V$ and $\dot{V}$.
\end{prop}
Each function $\omega_{i}$ depends on the previous gains $\left(\tilde{k}_{i-1},\cdots,\,\tilde{k}_{n}\right)$,
$\beta_{0,\,j}$, $\beta_{\infty,\,j}$, $d_{0}$, $d_{\infty}$,
$p_{0}$, $p_{\infty}$, $\kappa_{j}$, 
and $\theta_{j}$. Due to the recursive nature of the process, gains $\left(\tilde{k}_{j},\cdots,\,\tilde{k}_{n}\right)$ are appropriate
for the differentiator of order $n-j+1$.

\subsection{Convergence acceleration and scaling the Lipschitz constant $\Delta$}

\label{subsec:Convergence-acceleration}
Perform on system (\ref{eq:LevDif_n}), for arbitrary constants
$\alpha>0$ and $L>0$, the following scaling of the gains, 
\begin{equation}
\kappa_{i}\rightarrow\left(\frac{L^{n}}{\alpha}\right)^{\frac{d_{0}}{r_{0,\,i}}}\kappa_{i}\,,\,\theta_{i}\rightarrow\left(\frac{L^{n}}{\alpha}\right)^{\frac{d_{\infty}}{r_{\infty,\,i}}}\theta_{i}\,,\,k_{i}\rightarrow L^{i} k_{i}\,.\label{eq:GainScaling}
\end{equation}
It is easy to show that the linear state transformation
\[
e_{i}\rightarrow\frac{L^{n-i+1}}{\alpha}e_{i}\,,
\]
together with a time scaling $t\rightarrow Lt$, transforms the scaled
error system to system (\ref{eq:Dif_error}). This means that the convergence is accelerated and the Lipschitz constant increased as
\[
T\left( z_0 \right) \rightarrow \frac{1}{L}T\left( z_0 \right)\,,\quad  \Delta\rightarrow\alpha\Delta\,.
\]

Using the scaling \eqref{eq:GainScaling}, it is possible to assign an \emph{arbitrary} pair of 
(worst case) \emph{convergence time} $\bar{T}^{*}$ and Lipschitz constant $\Delta^{*}$ to the differentiator, following the
procedure: 

(i) Given $d_{0}<0<d_{\infty}$, 
$\kappa_{j}>0$ and $\theta_{j}>0$, fix a set of stabilizing gains
$k_{i}$ and the corresponding supported perturbation size $\Delta$,
using e.g. Proposition \ref{prop:GainCalculation}. 

(ii) Calculate
the corresponding fixed-convergence time $\bar{T}$, either by means
of (\ref{eq:FxTEstimation}) or by simulations. 

%For a desired convergence time $\bar{T}^{*}$ and Lipschitz constant $\Delta^{*}$
(iii) Select the scaling gains $\left(\alpha,\,L\right)$ of (\ref{eq:GainScaling})
as $\alpha\geq\nicefrac{\Delta^{*}}{\Delta}$ and $L\geq\nicefrac{\bar{T}^{*}}{\bar{T}}$. 

This procedure generalizes to an arbitrary order and arbitrary degrees
that proposed in \cite{SeeHai20} for the first order differentiator
with $d_{0}=-1$. Note that this scaling, using two parameters, is
novel also for the homogeneous case.

\subsection{Effect of noise and the perturbation $\delta\left(t\right)$}

In the presence of noise the estimation error cannot be zero asymptotically,
but it is uniformly and ultimately bounded. Moreover, when $d_{0}>-1$
and $\Delta>0$ the estimation error is also only uniformly and ultimately
bounded. This also happens when $d_{0}=-1$, $d_{\infty}>-1$ and
the differentiator gains are not sufficiently large to fully compensate
the effect of $\delta\left(t\right)$.
\begin{prop}
\label{prop:ISS}Let the hypothesis of Theorem \ref{thm:Differentiator}
be satisfied and select stabilizing gains $k_{i}$ for the differentiator
\eqref{eq:LevDif_n}. If $-1<d_{0}\leq d_{\infty}<\frac{1}{n-1}$
or $-1=d_{0}<d_{\infty}<\frac{1}{n-1}$ then the estimation error
system (\ref{eq:ErrorDyn}) is Input-to-State Stable (ISS), considering $\nu\left(t\right)$ and $\delta\left(t\right)$ as 
inputs.\qed
\end{prop}
It follows from Proposition \ref{prop:ISS} that if noise and perturbation
are bounded, then the estimation error $z$ will be also bounded,
and if $\left(\nu\left(t\right),\,\delta\left(t\right)\right)\rightarrow0$
then $e\left(t\right)\rightarrow0$. The precision for small noise
signals is determined by the $0$-limit approximation and it is therefore
identical to the one of the homogeneous differentiator of homogeneity
degree $d_{0}$ \cite{Lev03,Lev05,CruMor19,SanCru18}. In particular,
when $d_{0}=-1$, $\left|\nu\left(t\right)\right|\leq\epsilon$, $\left|\delta\left(t\right)\right|\leq\Delta$
the following inequalities are achieved in finite-time 
\[
\left|x_{i}\left(t\right)-f_{0}^{\left(i-1\right)}\left(t\right)\right|\leq\lambda_{i}\Delta^{\frac{i-1}{n}}\left|\epsilon\right|^{\frac{n-i+1}{n}}\,,\,\forall t\geq T\,.
\]

\section{Proof of the results: A Lyapunov approach}

\label{sec:LyapunovProof} We write (\ref{eq:ErrorDyn}) in compact
form as $\dot{z}\in F\left(z\right)+b\bar{\delta}\left(t\right)$,
where $b=\left[0,\cdots,\,0,\,1\right]^{T}\in\mathbb{R}^{n}$. Since
for $d_{0}=-1$ the function $\phi_{n}\left(z_{1}\right)$ is set-valued
due to the sign function, $F\left(z\right)$ is in general a set-valued
vector field which satisfies standard assumptions.

Since for $i=1,\cdots,\,n-1$, $r_{0,\,i}>0$ and $r_{\infty,\,i}>0$,
each function $\varphi_{i}\left(z_{i}\right)$ in (\ref{eq:DefVarphi_i})
is $\mathcal{C}$ on $\mathbb{R}$, $\mathcal{C}^{1}$ on $\mathbb{R}\setminus\left\{ 0\right\} $,
strictly increasing and surjective. Its inverse $\varphi_{i}^{-1}\left(z_{i}\right)$
is well-defined, $\mathcal{C}$ on $\mathbb{R}$, $\mathcal{C}^{1}$
on $\mathbb{R}\setminus\left\{ 0\right\} $, and also strictly increasing.
For $\varphi_{n}\left(z_{n}\right)$ the same is true if $d_{0}>-1$.
If $d_{0}=-1$ function $\varphi_{n}\left(z_{n}\right)=\kappa_{n}\left\lceil z_{n}\right\rfloor ^{0}+\theta_{n}\left\lceil z_{n}\right\rfloor ^{1+d_{\infty}}$
is discontinuous in $z_{n}=0$, and $\mathcal{C}^{1}$ on $\mathbb{R}\setminus\left\{ 0\right\} $. 

Since $d_{\infty}>d_{0}$, $\varphi_{i}\left(z_{i}\right)$ in (\ref{eq:DefVarphi_i})
is homogeneous in the $0$-limit and in the $\infty$-limit, with
approximating functions $\varphi_{i,\,0}\left(z_{i}\right)=\kappa_{i}\left\lceil z_{i}\right\rfloor ^{\frac{r_{0,\,i+1}}{r_{0,\,i}}}$
and $\varphi_{i,\,\infty}\left(z_{i}\right)=\theta_{i}\left\lceil z_{i}\right\rfloor ^{\frac{r_{\infty,\,i+1}}{r_{\infty,\,i}}}$,
respectively. For $i=1,\cdots,\,n-1$, the inverse $\varphi_{i}^{-1}\left(s\right)$
is also homogeneous in the $0$-limit and in the $\infty$-limit,
with approximating functions $\varphi_{i,\,0}^{-1}\left(s\right)=\frac{1}{\theta_{i}}\left\lceil s\right\rfloor ^{\frac{r_{\infty,\,i}}{r_{\infty,\,i+1}}}$
and $\varphi_{i,\,\infty}^{-1}\left(s\right)=\frac{1}{\kappa_{i}}\left\lceil s\right\rfloor ^{\frac{r_{0,\,i}}{r_{0,\,i+1}}}$,
respectively. Note also that when $d_{\infty}>0$, for $i=1,\cdots,\,n-1$,
$\varphi_{i,\,0}^{-1}\left(s\right)$ is homogeneous of negative degree,
and therefore $\varphi_{i}^{-1}\left(s\right)$ not differentiable
at $s=0$. However, $\left\lceil \varphi_{i}^{-1}\left(s\right)\right\rfloor ^{\mu}$
is differentiable at $s=0$ for every $\mu\geq\frac{r_{\infty,\,i}+d_{\infty}}{r_{\infty,\,i}}$,
and $\left|\varphi_{i}^{-1}\left(s\right)\right|^{\mu}$ for every
$\mu>\frac{r_{\infty,\,i}+d_{\infty}}{r_{\infty,\,i}}$.

For $i=1,\cdots,\,n-1$, functions $\phi_{i}$ in (\ref{eq:DefPhi_i}),
being compositions of $\varphi_{j}$, are $\mathcal{C}$ on $\mathbb{R}$,
$\mathcal{C}^{1}$ on $\mathbb{R}\setminus\left\{ 0\right\} $, strictly
increasing and surjective. In case $d_{0}=-1$, function $\phi_{n}\left(z_{1}\right)=\kappa_{n}\left\lceil z_{1}\right\rfloor ^{0}+\theta_{n}\left\lceil \varphi_{n-1}\circ\cdots\varphi_{2}\circ\varphi_{1}\left(z_{1}\right)\right\rfloor ^{\frac{r_{\infty,\,n+1}}{r_{\infty,\,n}}}$
is discontinuous. Since it is used in the Differential Inclusion (\ref{eq:ErrorDyn}),
using the Filippov's regularization procedure \cite{Fil88} it will
become an upper semi-continuous set-valued function, where the sign
function $\left\lceil s\right\rfloor ^{0}$ is defined as usually
for $s\neq0$, but for $s=0$ its values are an interval, i.e. $\left\lceil 0\right\rfloor ^{0}=\left[-1,\,1\right]\in\mathbb{R}$.
$\phi_{i}$'s are also homogeneous in the $0$-limit and in
the $\infty$-limit, with approximating functions 
\[
\phi_{i,0}\left(s \right) =  K_{i,0} \left\lceil s \right\rfloor ^{\frac{r_{0,i+1}}{r_{0,1}}} \,,
\phi_{i,\infty}\left( s \right)  = K_{i,\infty} \left\lceil s \right\rfloor ^{\frac{r_{\infty,i+1}}{r_{\infty,1}}},\,
\]
where $K_{i,0} = \prod_{j=1}^{i}\kappa_{j}^{\frac{r_{0,i+1}}{r_{0,j+1}}}$, $K_{i,\infty} = \prod_{j=1}^{i}\theta_{j}^{\frac{r_{\infty,i+1}}{r_{\infty,j+1}}}$. %
For $\delta\left(t\right)\equiv0$ system (\ref{eq:ErrorDyn}) is
bl-homogeneous with homogeneity degrees $d_{0}$ and $d_{\infty}$
and weights $\mathbf{r}_{0}=\left[r_{0,\,1},\cdots,\,r_{0,\,n}\right]$
and $\mathbf{r}_{\infty}=\left[r_{\infty,\,1},\cdots,\,r_{\infty,\,n}\right]$
as in (\ref{eq:Def_r}).

Since $\varphi_{n}\left(z_{n}\right)$ is not involved in the definition
of $Z_{i}$, it has to satisfy weakened conditions compared to the
other functions $\varphi_{i}$. From the properties of functions $\varphi_{i}$
it follows that $Z_{i}$ is $\mathcal{C}$ on $\mathbb{R}$. For $Z_{i}$
to be $\mathcal{C}^{1}$ on $\mathbb{R}$ the powers in (\ref{eq:DefZi})
have to be sufficiently large, what is the case if (\ref{eq:Cond_p1})
is fulfilled. Note that if (\ref{eq:Cond_p2}) is met, $Z_{i}$ is
also bl-homogeneous with approximations $Z_{i,\,0}\left(z_{i},\,z_{i+1}\right)$, given by the first term in \eqref{eq:DefZi} with $\xi = \frac{1}{\theta_{i}}\left\lceil z_{i+1}\right\rfloor ^{\frac{r_{\infty,\,i}}{r_{\infty,\,i+1}}}$; and  $Z_{i,\,\infty}\left(z_{i},\,z_{i+1}\right)$, given by the second term in \eqref{eq:DefZi} with $\xi = \frac{1}{\kappa_{i}}\left\lceil z_{i+1}\right\rfloor ^{\frac{r_{0,\,i}}{r_{0,\,i+1}}}$. %
Moreover, functions $Z_{i}\left(z_{i},\,z_{i+1}\right)$ are nonnegative.
\begin{lem}
$Z_{i}\left(z_{i},\,z_{i+1}\right)\geq0$ for every $i=1,\cdots,\,n$
and $Z_{i}\left(z_{i},\,z_{i+1}\right)=0$ if and only if $\varphi_{i}\left(z_{i}\right)=z_{i+1}$.
\end{lem}
\begin{IEEEproof}
From Young's inequality it follows that
\[
z_{i}\left\lceil \xi \right\rfloor ^{\frac{p_{0}-r_{0,\,i}}{r_{0,\,i}}} \leq \frac{r_{0,\,i}}{p_{0}}\left|z_{i}\right|^{\frac{p_{0}}{r_{0,\,i}}}+\left(1-\frac{r_{0,\,i}}{p_{0}}\right)\left| \xi \right|^{\frac{p_{0}}{r_{0,\,i}}}\,.
\]
From this and (\ref{eq:DefZi}) it follows that $Z_{i}\left(z_{i},\,z_{i+1}\right)\geq0$.
\end{IEEEproof}
The partial derivatives of $Z_{i}\left(z_{i},\,z_{i+1}\right)$, for
which we introduce the symbols $\sigma_{i}$ and $s_{i}$, are given
by 
\begin{gather}
\sigma_{i}\left(z_{i},\,z_{i+1}\right) \triangleq \frac{\partial Z_{i}\left(z_{i},\,z_{i+1}\right)}{\partial z_{i}} =  \qquad \qquad \qquad  \label{eq:DefSigma_i} \\
\sum_{j=\{0,\,\infty \}} \beta_{j,i}\left(\left\lceil z_{i}\right\rfloor ^{\frac{p_{j}-r_{j,i}}{r_{j,i}}}-\left\lceil \varphi_{i}^{-1}\left(z_{i+1}\right)\right\rfloor ^{\frac{p_{j}-r_{j,\,i}}{r_{j,\,i}}}\right) \nonumber  \\
%\end{gather}
%\begin{gather}
s_{i}\left(z_{i},\,z_{i+1}\right) \triangleq \frac{\partial Z_{i}\left(z_{i},\,z_{i+1}\right)}{\partial z_{i+1}} =  \qquad \qquad \qquad   \label{eq:DefS_i} \\
\sum_{j=\{0,\,\infty \}} -\beta_{j,i}\frac{p_{j}-r_{j,i}}{r_{j,i}}\left( z_{i} - \xi_{i} \right) \left| \xi_{i} \right|^{\frac{p_{j}-2r_{j,i}}{r_{j,i}}}\frac{\partial \xi_{i} }{\partial z_{i+1}} \nonumber 
\end{gather}
where $\xi_{i} = \varphi_{i}^{-1}\left(z_{i+1}\right)$. %
Note that $s_{n}\left(z_{n},\,z_{n+1}\right)\equiv0$, and that functions
$\sigma_{i}\left(z_{i},\,z_{i+1}\right)$ and $s_{i}\left(z_{i},\,z_{i+1}\right)$
are $\mathcal{C}$ on $\mathbb{R}$, bl-homogeneous of degrees $p_{0}-r_{0,\,i}$,
$p_{0}-r_{0,\,i+1}$ for the $0$-approximation and $p_{\infty}-r_{\infty,\,i}$,
$p_{\infty}-r_{\infty,\,i+1}$ for the $\infty$-approximation, respectively.
Futhermore, for $i=1,\cdots,\,n-1$, $\{\sigma_{i}=0 \}= \{s_{i}=0 \}$, i.e. $\sigma_i$, $s_i$ are zero where $Z_{i}$
achieves its minimum $Z_{i}=0$.

$V$ is bl-homogeneous of degrees $p_{0}$ and $p_{\infty}$ and $\mathcal{C}^{1}$ on $\mathbb{R}$. It is also non negative, since it is a positive
combination of non negative terms. Moreover, $V$ is
positive definite since $V\left(z\right)=0$ only if all $Z_{i}=0$,
what only happens at $z=0$. Due to bl-homogeneity it is radially
unbounded \cite{Bhat2005}.

For calculation of the time derivative of $V$ along the trajectories
of (\ref{eq:ErrorDyn}), we consider first the nominal situation in
which $\delta\left(t\right)\equiv0$ if $d_{0}\neq-1$ and $\left|\delta\left(t\right)\right|\leq\Delta$
when $d_{0}=-1$. In that case
\begin{equation}
\dot{V}\left(z\right)  \in  W\left(z\right)\,,\label{eq:Vdot}
\end{equation}
where
\begin{align}
W\left(z\right) & =  -\tilde{k}_{1}\sigma_{1} \left(\phi_{1}\left(z_{1}\right)-z_{2}\right)\nonumber \\
 &   -\sum_{j=2}^{n-1}\tilde{k}_{j}\left[s_{j-1}+\sigma_{j} \right]\left(\phi_{j}\left(z_{1}\right)-z_{j+1}\right)\label{eq:W}\\
 &   -\tilde{k}_{n}\left[s_{n-1}+\sigma_{n} \right]\left(\phi_{n}\left(z_{1}\right)-\frac{\Delta}{k_{n}}\left[-1,\,1\right]\right)\nonumber 
\end{align}
where $\left[-1,\,1\right]\in\mathbb{R}$ is an interval and we omit the variable dependence of $\sigma_i$ and $s_i$. Since the
set-valued vector field $F\left(z\right)$ on the right-hand side
of (\ref{eq:ErrorDyn}) with $\bar{\delta}\left(t\right)\equiv0$
is upper semi-continuous and the gradient $\nabla V\left(z\right)$
of $V\left(z\right)$ is continuous, by \cite[Lemma 6]{CruMor20b}
the set-valued function $W\left(z\right)=\nabla V\left(z\right)F\left(z\right)$
is also upper semi-continuous, and the single-valued function $W^{*}\left(z\right)=\max\left\{ W\left(z\right)\right\} $
is upper semi-continuous and bl-homogeneous of degrees $p_{0}+d_{0}$
and $p_{\infty}+d_{\infty}$, respectively. Moreover, $W\left(0\right)=0$
since $\nabla V\left(0\right)=0$. We want to show that there exist
values of $\tilde{k}_{i}>0$ such that $W\left(z\right)<0$, i.e.
$W$ is negative definite. This is equivalent to showing that $W^{*}\left(z\right)<0$.
We note that when $-1<d_{0}$ function $W$ is indeed single-valued
and continuous, and thus $W^{*}\left(z\right)=W\left(z\right)$ (recall
that in this case $\Delta=0$). When $d_{0}=-1$, $W$ is set-valued
because the term 
$\phi_{n}\left(z_{1}\right)-\frac{\Delta}{k_{n}}\left[-1,\,1\right]=\kappa_{n}\left\lceil z_{1}\right\rfloor ^{0}+\theta_{n}\left\lceil \phi_{n-1}\left(z_{1}\right)\right\rfloor ^{\frac{r_{\infty,\,n+1}}{r_{\infty,\,n}}}-\frac{\Delta}{k_{n}}\left[-1,\,1\right]$
is set-valued. To simplify the development in what follows, we will
write simply $\phi_{n}\left(z_{1}\right)$ but it is meant $\phi_{n}\left(z_{1}\right)-\frac{\Delta}{k_{n}}\left[-1,\,1\right]$.

For the proof we will use the following property of upper semi-continuous,
bl-homogeneous single-valued functions, proven in \cite{CruMor20b}.
\begin{lem}
\label{lemMorEmus} Let $\eta:\mathbb{R}^{n}\rightarrow\mathbb{R}$
and $\gamma:\mathbb{R}^{n}\rightarrow\mathbb{R}_{\leq0}$ be two upper
semicontinuous (u.s.c.) single-valued \emph{bl}-homogeneous functions,
with the same weights $\mathbf{r}_{0}$ and $\mathbf{r}_{\infty}$,
degrees $m_{0}$ and $m_{\infty}$, and approximating functions $\eta_{0}$,
$\eta_{\infty}$ and $\gamma_{0}$, $\gamma_{\infty}$, which are u.s.c. Suppose
that $\forall x\in\mathbb{R}^{n}$, 
$\gamma\left(x\right)\leq0$, $\gamma_{0}\left(x\right)\leq0$, $\gamma_{\infty}\left(x\right)\leq0$. 
If
$\gamma\left(x\right)=0\wedge x\neq0  \Rightarrow\eta\left(x\right)<0$, 
$\gamma_{0}\left(x\right)=0\wedge x\neq0 \Rightarrow\eta_{0}\left(x\right)<0$, 
$\gamma_{\infty}\left(x\right)=0\wedge x\neq0 \Rightarrow\eta_{\infty}\left(x\right)<0$ %
then there are constants $\lambda^{\ast}\in\mathbb{R}$, $c_{0}>0$,
and $c_{\infty}>0$ such that for all $\lambda\geq\max\{\lambda_{0},\lambda_{\infty}\}$,
$\lambda_{0}\geq\lambda^{\ast}$, $\lambda_{\infty}\geq\lambda^{\ast}$
and for all $x\in\mathbb{R}^{n}\setminus\{0\}$, 
\begin{align*}
\eta\left(x\right)+\lambda\gamma\left(x\right) & \leq  -c_{0}\left\Vert x\right\Vert _{\mathbf{r}_{0},\,p}^{m_{0}}-c_{\infty}\left\Vert x\right\Vert _{\mathbf{r}_{\infty},\,p}^{m_{\infty}},\\
\eta_{0}\left(x\right)+\lambda\gamma_{0}\left(x\right) & \leq  -c_{0}\left\Vert x\right\Vert _{\mathbf{r}_{0},\,p}^{m_{0}},\\
\eta_{\infty}\left(x\right)+\lambda\gamma_{\infty}\left(x\right) & \leq  -c_{\infty}\left\Vert x\right\Vert _{\mathbf{r}_{\infty},\,p}^{m_{\infty}}\,. \qquad \qquad \qed
\end{align*}
%\qed 
\end{lem}
To show that $W\left(z\right)<0$ we exploit its structure. So consider
the values of $W$ restricted to some hypersurfaces: for $i=1,\cdots,\,n-1$
\[
\mathcal{Z}_{1}  =  \left\{ \varphi_{1}\left(z_{1}\right)=z_{2}\right\}  \cdots 
\mathcal{Z}_{i}  =  \mathcal{Z}_{i-1}\cap\left\{ \varphi_{i}\left(z_{i}\right)=z_{i+1}\right\} .
\]
These sets are clearly related as $\mathcal{Z}_{n-1}\subset \cdots  \subset\mathcal{Z}_{i} \subset\cdots\subset\mathcal{Z}_{1}\subset\mathbb{R}^{n}$.
Note that on $\mathcal{Z}_{i}$ functions $\sigma_{i}$ and $s_{i}$
vanish, i.e. $\sigma_{i}=s_{i}=0$, and therefore they also vanish
on $\mathcal{Z}_{j}$, for every $j>i$. Let $W_{i}=W_{\mathcal{Z}_{i}}$
represent the value of $W\left(z\right)$ restricted to the the manifold
$\mathcal{Z}_{i}$. We can obtain the value of $W_{1}$ by replacing
in $W\left(z\right)$ the variable $z_{1}$ by $z_{1}=\varphi_{1}^{-1}\left(z_{2}\right)$,
so that $W_{1}$ becomes a function of $\left(z_{2},\cdots,\,z_{n}\right)$.
In general, we obtain the value of $W_{i}$ , for $i=1,\dots,\,n-1$,
by replacing in $W\left(z\right)$ the variables $\left(z_{1},\cdots,\,z_{i}\right)$
by its values in terms of $z_{i+1}$, so that $W_{i}$ becomes a function
of $\bar{z}_{i+1}\triangleq\left(z_{i+1},\cdots,\,z_{n}\right)$.
%For example, the value of the expression $\left(\phi_{i}\left(z_{1}\right)-z_{i+1}\right)$ on $\mathcal{Z}_{i-1}$ is $\left(\varphi_{i}\left(z_{i}\right)-z_{i+1}\right)$.

The first term in $W\left(z\right)$ (\ref{eq:W}) is non positive,
i.e., $\sigma_{1}\left(z_{1},\,z_{2}\right)\left(\phi_{1}\left(z_{1}\right)-z_{2}\right)\leq0$,
and it vanishes on $\mathcal{Z}_{1}$. Evaluating $W\left(z\right)$
on $\mathcal{Z}_{1}$ we obtain (recall that $s_{1}=0$, $\phi_{i}\left(z_{1}\right)=\varphi_{i}\circ\cdots\varphi_{2}\circ\varphi_{1}\left(z_{1}\right)$
and $z_{n+1}\equiv0$)
\begin{multline*}
W_{1}\left(\bar{z}_{2}\right)  =  -\tilde{k}_{2}\sigma_{2}\left(z_{2},\,z_{3}\right)\left(\varphi_{2}\left(z_{2}\right)-z_{3}\right) \\
  -\sum_{j=3}^{n}\tilde{k}_{j}\left[s_{j-1}+\sigma_{j} \right]\left(\varphi_{j}\circ\cdots\circ\varphi_{2}\left(z_{2}\right)-z_{j+1}\right)\,.
\end{multline*}
Note that $W_{1}\left(\bar{z}_{2}\right)$ has the same structure
as $W\left(z\right)$. Its first term is non positive, i.e., $\sigma_{2}\left(z_{2},\,z_{3}\right)\left(\varphi_{2}\left(z_{2}\right)-z_{3}\right)\leq0$,
and it vanishes on $\mathcal{Z}_{2}$. Evaluating $W_{1}\left(\bar{z}_{2}\right)$
on $\mathcal{Z}_{2}$ we obtain $W_{2}\left(\bar{z}_{3}\right)$.
Applying this procedure recursively, and using the facts that $s_{i}=0$
on $\mathcal{Z}_{i}$ and $\phi_{i}\left(z_{1}\right)=\varphi_{i}\circ\cdots\varphi_{2}\circ\varphi_{1}\left(z_{1}\right)$,
we find that, for $i=1,\dots,\,n-1$,
\begin{multline*}
W_{i}\left(\bar{z}_{i+1}\right) = -\tilde{k}_{i+1}\sigma_{i+1}\left(\varphi_{i+1}\left(z_{i+1}\right)-z_{i+2}\right)\\
  -\sum_{j=i+2}^{n}\tilde{k}_{j}\left[s_{j-1}+\sigma_{j} \right]\left(\varphi_{j}\circ\cdots\circ\varphi_{i+1}\left(z_{i+1}\right)-z_{j+1}\right)\,.
\end{multline*}
Note that the first term of $W_{i}\left(\bar{z}_{i+1}\right)$ is
non positive, i.e., $\sigma_{i+1}\left(z_{i+1},\,z_{i+2}\right)\left(\varphi_{i+1}\left(z_{i+1}\right)-z_{i+2}\right)\leq0$,
and it vanishes on $\mathcal{Z}_{i+1}$. For $i=n-1$ the value of
$W_{n-1}\left(\bar{z}_{n}\right)$ is given by (recall that $z_{n+1}\equiv0$)
\begin{multline*}
W_{n-1}\left(\bar{z}_{n}\right)  % = -\tilde{k}_{n}\sigma_{n}\left(z_{n},\,z_{n+1}\right)\varphi_{n}\left(z_{n}\right)\\
 = -\tilde{k}_{n}\left(\beta_{0,n}\left\lceil z_{n}\right\rfloor ^{p_{0}-1}+\beta_{\infty,n}\left\lceil z_{n}\right\rfloor ^{p_{\infty}-1}\right)  \\ \times \left(\kappa_{n}\left\lceil z_{n}\right\rfloor ^{1+d_{0}}+\theta_{n}\left\lceil z_{n}\right\rfloor ^{1+d_{\infty}}-\frac{\Delta}{k_{n}}\left[-1,\,1\right]\right)\,,
\end{multline*}
where we have used (\ref{eq:DefSigma_i}) and (\ref{eq:DefVarphi_i}). %
Here we distinguish two cases: (i) $-1<d_{0}$ and (ii) $d_{0}=-1$.

When $-1<d_{0}$, $\Delta=0$, $W$ is single-valued and continuous,
$W_{n-1}\left(\bar{z}_{n}\right)$ is bl-homogeneous and it is negative
for any $\tilde{k}_{n}>0$. 

In case $d_{0}=-1$, if $\tilde{\Delta} \triangleq \frac{\Delta}{\kappa_{n}k_{n}}<1$ the following
equality is satisfied for the set-valued map 
\[
\kappa_{n}\left( \left\lceil z_{n}\right\rfloor ^{0}- \tilde{\Delta} \left[-1,\,1\right] \right) = \kappa_{n} \left\lceil z_{n}\right\rfloor ^{0}\left[1- \tilde{\Delta},\,1+\tilde{\Delta}\right]\,,
\]
so that it is clear that for any $\nu>0$, $\kappa_{n}\left(\left\lceil z_{n}\right\rfloor ^{0}- \tilde{\Delta} \left[-1,\,1\right]\right)\left\lceil z_{n}\right\rfloor ^{\nu}>0$
for $z_{n}\neq0$ and $\kappa_{n}\left(\left\lceil z_{n}\right\rfloor ^{0}- \tilde{\Delta} \left[-1,\,1\right]\right)\left\lceil z_{n}\right\rfloor ^{\nu}=0$
for $z_{n}=0$. Therefore, with $W_{n-1}^{*}\left(\bar{z}_{n}\right) = \max\left\{ W_{n-1}\left(\bar{z}_{n}\right)\right\}$,
\begin{multline*}
W_{n-1}^{*}\left(\bar{z}_{n}\right) = -\tilde{k}_{n}\left(\beta_{0,n}\left|z_{n}\right|^{p_{0}-1}+\beta_{\infty,n}\left|z_{n}\right|^{p_{\infty}-1}\right) \\
\times \left(\kappa_{n}\left(1-\frac{\Delta}{\kappa_{n}k_{n}}\right)+\theta_{n}\left|z_{n}\right|^{1+d_{\infty}}\right)\,.
\end{multline*}
Function $W_{n-1}^{*}$ is
single-valued, upper semi-continuous, bl-homogeneous and negative
definite for any $\tilde{k}_{n}>0$.

In all cases, the same is true for its homogeneous approximations (as shown in \cite{CruMor16,SanCru18}).
$W_{n-2}\left(\bar{z}_{n-1}\right)$ is also bl-homogeneous. According
to Lemma \ref{lemMorEmus} we conclude that $W_{n-2}\left(\bar{z}_{n-1}\right)$
can be rendered negative definite (in $\mathcal{Z}_{n-2}$) by selecting
$\tilde{k}_{n-1}>0$ sufficiently large. Since $W_{i}\left(\bar{z}_{i+1}\right)$
is bl-homogeneous and the conditions of Lemma \ref{lemMorEmus} are
satisfied for $W_{i}\left(\bar{z}_{i+1}\right)$ and its homogeneous
approximations (as shown in \cite{CruMor16,SanCru18}),
we conclude that $W_{i-1}\left(\bar{z}_{i}\right)$ can be rendered
negative definite (in $\mathcal{Z}_{i-1}$) selecting $\tilde{k}_{i}>0$
sufficiently large. %
Applying the argument recursively, we conclude that there exist positive
values of $\left(\tilde{k}_{1},\cdots,\,\tilde{k}_{n}\right)$ such
that $W\left(z\right)<0$.

Moreover, using \cite[Corollary 2.15]{AndPra08} (see also \cite[Lemma 10]{CruMor20b})
the inequality (\ref{eq:LyapIneq}) follows. Using inequality
(\ref{eq:LyapIneq}) satisfied by the Lyapunov function we obtain,
as a direct consequence of \cite[Lemma 3]{CruMor20b}, the estimation
of the convergence time given by \eqref{eq:FxTEstimation}. 

\subsection{Gain calculation}

\label{subsec:Gain-calculation}The gains $\tilde{k}_{i}$ are calculated
backwards from $i=n,\,n-1,\,\cdots,\,2,\,1$ such that $W_{i}\left(\bar{z}_{i+1}\right)>0$.
This will be the case if they are chosen as given in Proposition \ref{prop:GainCalculation},
where functions $\omega_{i}$ are defined as (recall that $\varphi_{n}\left(s\right)$
should be replaced by $\varphi_{n}\left(s\right)-\frac{\Delta}{k_{n}}\left[-1,\,1\right]$)
: 
\begin{gather}
\omega_{i}\left(\tilde{k}_{i+1},\cdots,\,\tilde{k}_{n}\right)\triangleq  \max_{\left(z_{i},\cdots,\,z_{n}\right)\in\mathbb{R}^{n-i+1}} \left\{ \right. \qquad  \label{Def:Omegai}\\
\left. \frac{\sum_{j=i+1}^{n}\tilde{k}_{j}\left[s_{j-1}+\sigma_{j}\right]\left(z_{j+1}-\varphi_{j}\circ\cdots\circ\varphi_{i}\left(z_{i}\right)\right)}{\sigma_{i}\left(z_{i},\,z_{i+1}\right)\left(\varphi_{i}\left(z_{i}\right)-z_{i+1}\right)} \right\}. \nonumber
\end{gather}
These maximizations are well-posed, as it is shown in the previous
steps of the proof. %For $\tilde{k}_{i}$ the maximization has to be performed on the set $\left(z_{i},\cdots,\,z_{n}\right)\in\mathbb{R}^{n-i+1}$. 

\subsection{ISS and effect of noise}

We prove the Proposition \ref{prop:ISS} using function $V\left(z\right)$ (\ref{eq:DefLF}), although it 
is not a Lyapunov Function for all possible stabilizing gains $k_{i}$. However, 
the converse Lyapunov theorem for bl-homogeneous differential inclusions
\cite[Theorem 1]{CruMor20b} assure the existence of an appropriate
one. For it, the calculations are similar to those performed with (\ref{eq:DefLF}).
If we consider the effect of noise, in the continuous case, i.e. $d_{0}>-1$,
$\dot{V}$ can be written as
\begin{gather*}
\dot{V}\left(z\right) = \frac{1}{2}W\left(z\right)+R\left(z,\,\nu,\,\delta\right)\\
R\left( \cdot \right)  \triangleq  \frac{1}{2}W -\sum_{j=1}^{n}\tilde{k}_{j}\left[s_{j-1}+\sigma_{j} \right]\left(\phi_{j}\left(z_{1}-\nu\right)-\phi_{j}\left(z_{1}\right)\right)\\
 +\tilde{k}_{n}\left[s_{n-1}\left(z_{n-1},\,z_{n}\right)+\sigma_{n}\left(z_{n}\right)\right]\bar{\delta}\left(t\right) \,.
\end{gather*}
Define $Q\left(\nu,\,\bar{\delta}\right) = \left|\nu\right|^{\frac{p_{0}+d_{0}}{r_{0,\,1}}}+\left|\nu\right|^{\frac{p_{\infty}+d_{\infty}}{r_{\infty,\,1}}}+\left|\bar{\delta}\right|^{\frac{p_{0}+d_{0}}{1+d_{0}}}+\left|\bar{\delta}\right|^{\frac{p_{\infty}+d_{\infty}}{1+d_{\infty}}}$.  $W$, $R$ and $Q$ are continuous bl-homogeneous functions of degrees $p_{0}+d_{0}$
and $p_{\infty}+d_{\infty}$, and $Q$ is non negative. Furthermore, the function
$R\left(z,\,\nu,\,\bar{\delta}\right)-\gamma Q\left( \nu,\,\bar{\delta} \right)$ is also continuous and bl-homogeneous and $R\left(z,\,0,\,0\right)<0$
for $z\neq0$. The same is true for the homogeneous approximations.
And therefore, using \cite[Corollary 2.15]{AndPra08} (see also \cite[Lemma 10]{CruMor20b})
we conclude that there exists $\gamma>0$ such that %
$R\left(z,\,\nu,\,\bar{\delta}\right)\leq \gamma Q\left( \nu,\,\bar{\delta} \right)$. 
And thus
$\dot{V}\left(z\right)  \leq  -\frac{1}{2}\eta_{0}V^{\frac{p_{0}+d_{0}}{p_{0}}}\left(z\right)-\frac{1}{2}\eta_{\infty}V^{\frac{p_{\infty}+d_{\infty}}{p_{\infty}}}\left(z\right)+\gamma Q\left( \nu,\,\bar{\delta} \right)$, 
implying ISS by standard arguments. For the discontinuous case, $d_0=-1$, we can use the procedure used in \cite{SanCru18}.

\section{Example}

\label{sec:Examples} We perform some simulations with the bl-homogeneous second order differentiator ($n=3$), with $d_{0}=-1$ and  $d_{\infty}=\frac{1}{5}$. %
The signal to be differentiated is $f_{0}\left(t\right)=\frac{1}{2}\sin\left(\frac{1}{2}t\right)+\frac{1}{2}\cos\left(t\right)$,
for which $\Delta=\frac{5}{8}$, the internal gains $\kappa_{i}=\theta_{i}=1$
for $i=1,\,2,\,3$ and the gains $k_{1}=3$, $k_{2}=1.5\sqrt{3}$,
$k_{3}=1.1$. Two values of the scaling parameters were selected $\left(\alpha,\,L\right)=\left(1,\,1\right)$
and $\left(\alpha,\,L\right)=\left(1,\,2\right)$.

\begin{figure*}
     \centering
     \begin{subfigure}[b]{0.5\textwidth}
         \centering
         \includegraphics[width=\textwidth]{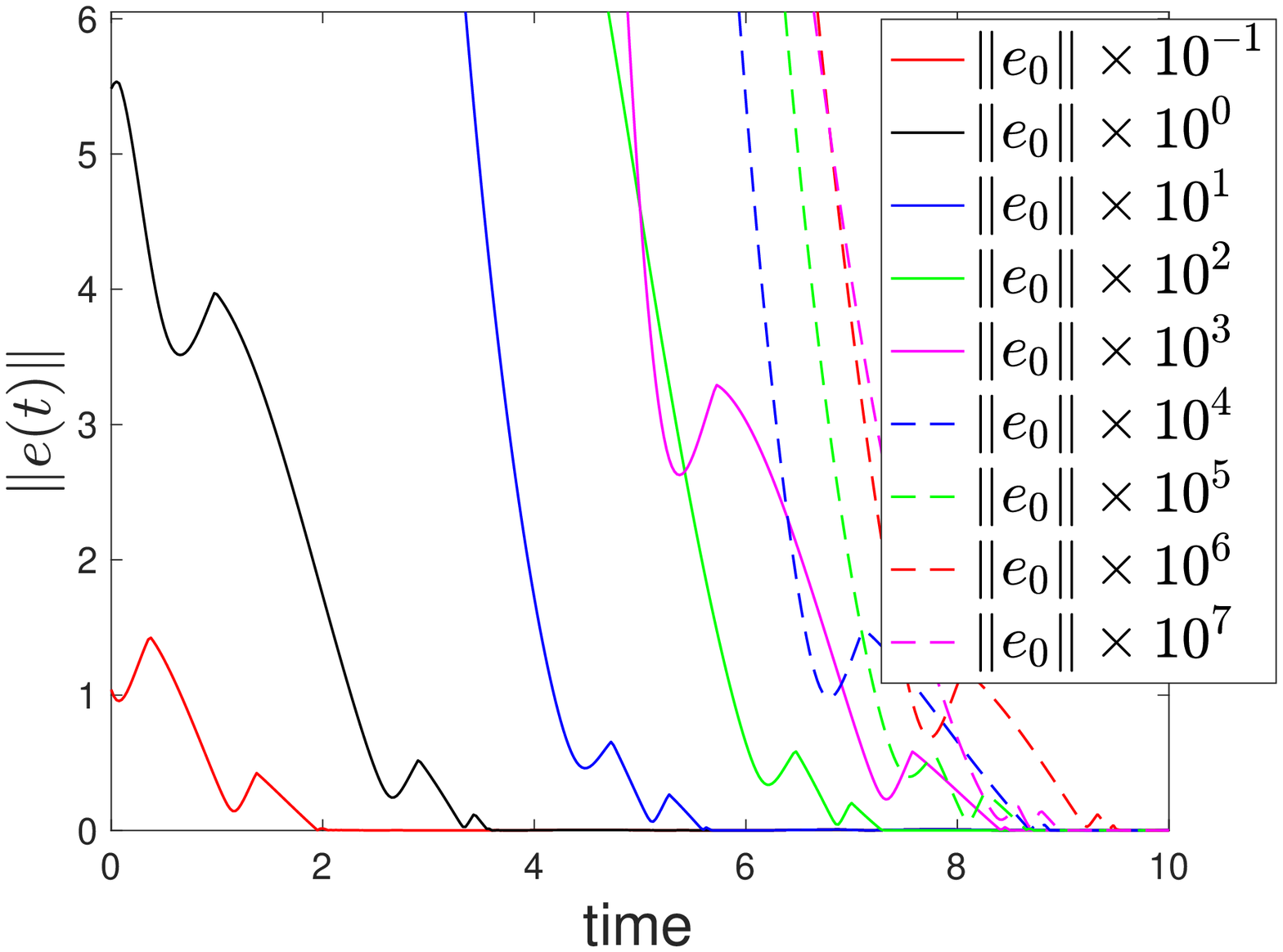}
         \caption{$\left\Vert e\left(t\right)\right\Vert $ for different
$\left\Vert e_0 \right\Vert $ in a logarithmic succession and  $L=1$.}
         \label{Fig1}
     \end{subfigure}
     \hfill
     \begin{subfigure}[b]{0.5\textwidth}
         \centering
         \includegraphics[width=\textwidth]{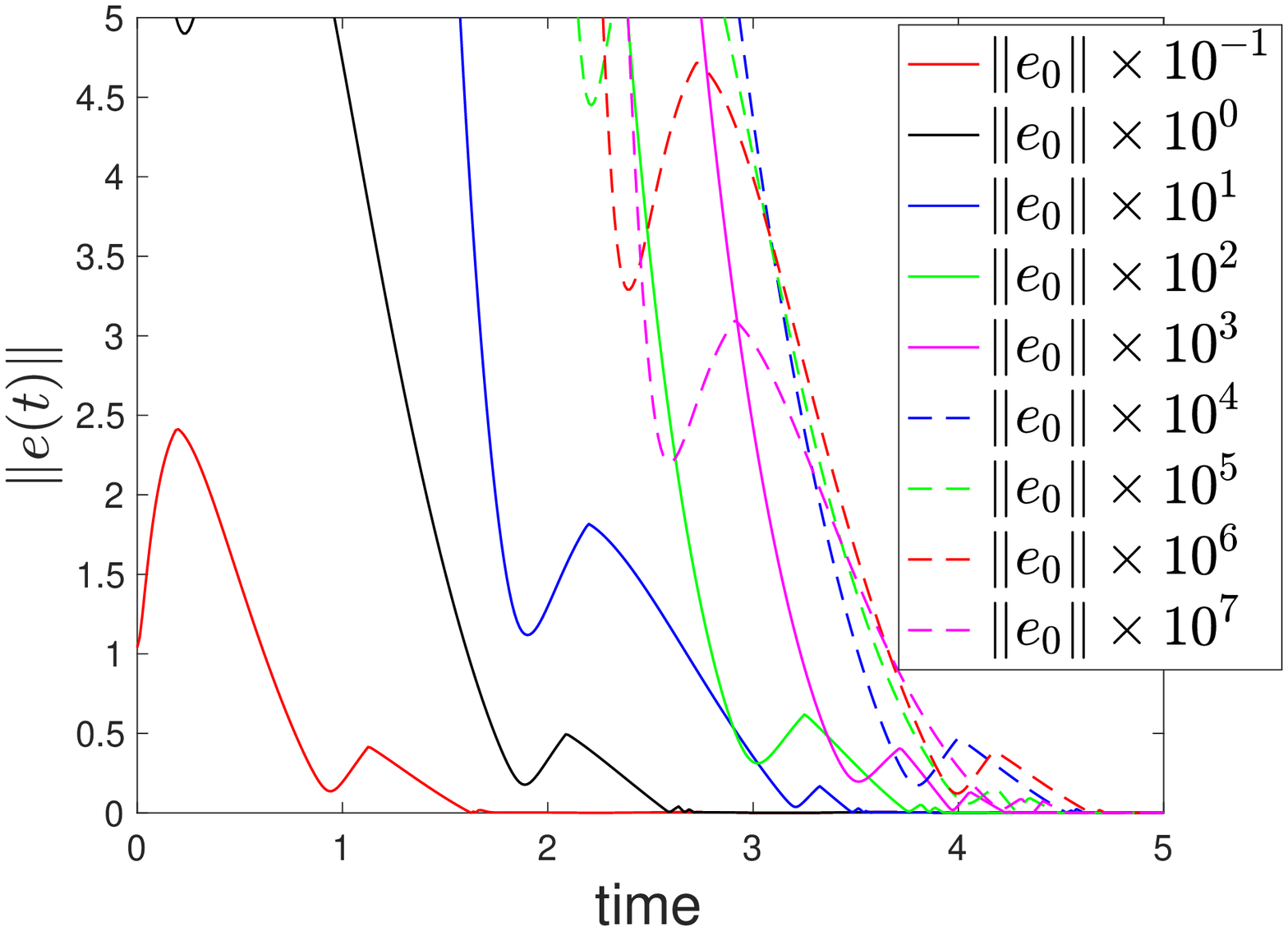}
         \caption{$\left\Vert e\left(t\right)\right\Vert $ for different
$\left\Vert e_0 \right\Vert $ in a logarithmic succession and  $L=2$.}
         \label{Fig2}
     \end{subfigure}
     \hfill
     \begin{subfigure}[b]{0.5\textwidth}
         \centering
         \includegraphics[width=\textwidth]{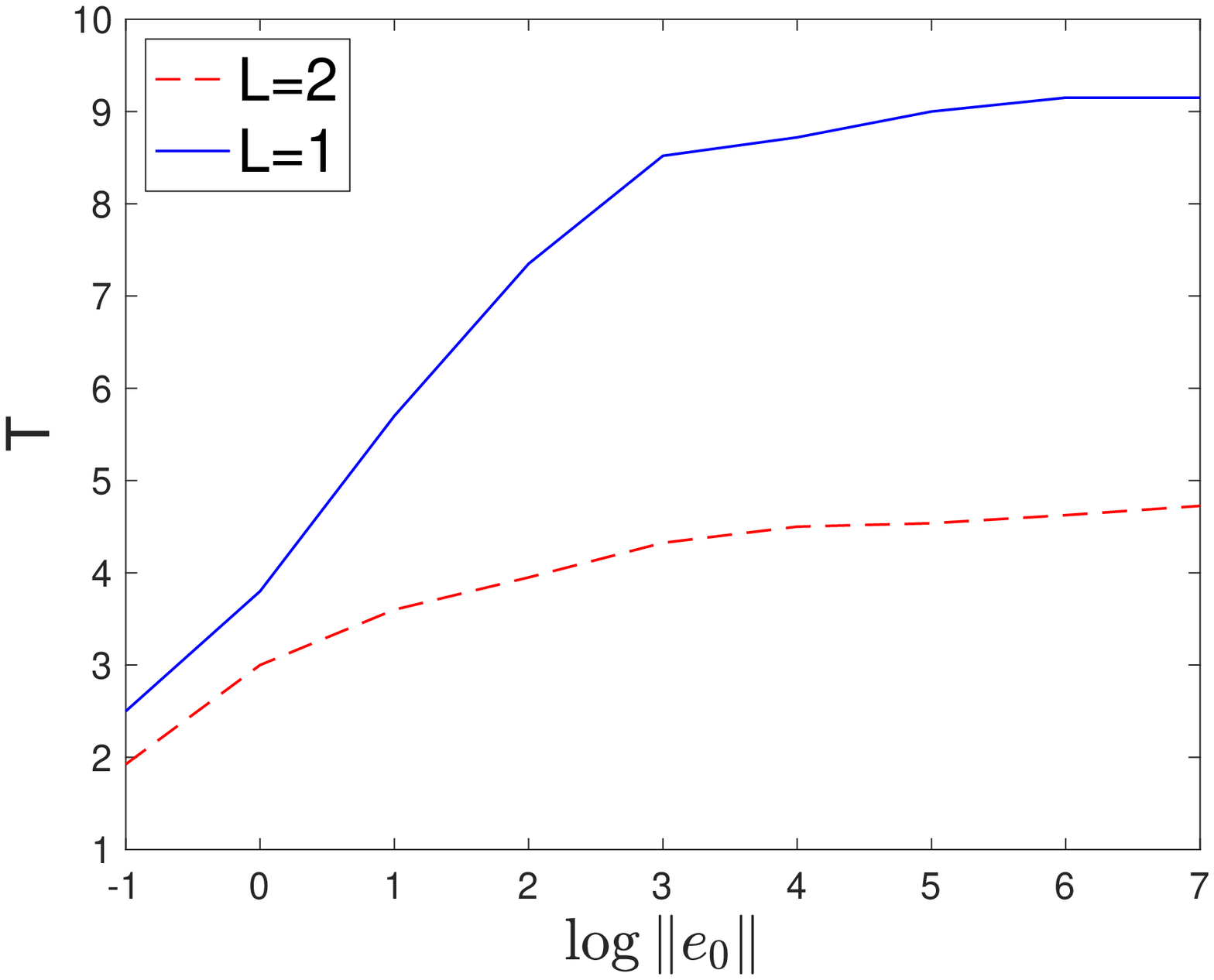}
         \caption{Convergence time versus the logarithmus of the initial condition $\left\Vert e_{0}\right\Vert $.}
         \label{Fig3}
     \end{subfigure}
        \caption{Time behavior of the estimation error norm $\left\Vert e\left(t\right)\right\Vert $  and Convergence time.}
%        \label{Fig3}
\end{figure*}

In Figures \ref{Fig1} and \ref{Fig2} the norm $\left\Vert e\left(t\right)\right\Vert $
of the estimation error is presented for different initial conditions
$e_{0}=\left[1,\,-5,\,1\right]\times10^{p}$, for $p=-1,\,0,\,1,\cdots,\,7$,
with $L=1$ in Figure \ref{Fig1} and $L=2$ in Figure \ref{Fig2}.
It is apparent from these graphs, that despite of a change in the
initial conditions of $8$ orders of magnitude, the convergence time
does not increase accrodingly, and it approaches an asymptote. This
is illustrated in Figure \ref{Fig3}, where the convergence time versus
the (logaritmic) value of the initial condition is shown. The figures
also show that by doubling the scaling parameter $L$ from $L=1$
to $L=2$, the convergence time is halved. In fact, using the parameter
$L$ any arbitrary convergence time can be attained.

\section{Conclusions}

\label{sec:Conclusions} We have proposed Fixed-Time converging exact and robust differentiators. In particular, Levant's discontinuous differentiator is extended with higher order terms, so that its convergence time is independent of the initial estimation error and can be arbitrarily assigned. Moreover, a full family of continuous differentiators are also studied, in a unified Lyapunov framework. We use the concept of homogeneity in the bi-limit, developed in \cite{AndPra08}, and the recursive observer design, for our objective. 

%\section*{Bibliography}

%\bibliographystyle{IEEEtran}
%\bibliography{/Users/jmorenop/Dropbox/DocumentosOld/ArticulosEnProceso/Biblio/MisArticulos,/Users/jmorenop/Dropbox/DocumentosOld/ArticulosEnProceso/Biblio/biblioMot,/Users/jmorenop/Dropbox/DocumentosOld/ArticulosEnProceso/Biblio/biblioam}

\end{document}